\documentclass[11pt]{amsart}
\usepackage{amsmath,amssymb,latexsym,soul,cite,amsthm,color,enumitem,graphicx,tikz, mathtools,microtype}
\usepackage{float}

\usepackage[colorlinks=true, urlcolor=blue, citecolor=red, linkcolor=blue]{hyperref}
\usepackage[english]{babel}
\usepackage[left=2.5cm,right=2.5cm,top=2.7cm,bottom=2.7cm]{geometry}

\numberwithin{equation}{section}

\newtheorem{theorem}{Theorem}[section]
\theoremstyle{plain}

\theoremstyle{plain}
\newtheorem{proposition}[theorem]{Proposition}
\theoremstyle{plain}
\newtheorem{corollary}[theorem]{Corollary}

\theoremstyle{definition}
\newtheorem{remark}[theorem]{Remark}

\newcommand{\N}{{\mathbb N}}

\newcommand{\R}{{\mathbb R}}

\newcommand{\beq}{\begin{equation}}
\newcommand{\eeq}{\end{equation}}
\renewcommand{\le}{\leqslant}
\renewcommand{\ge}{\geqslant}

\newcommand{\supp}{\operatorname{supp}}

\newcommand{\osc}{\operatornamewithlimits{osc}}

\newcommand{\dist}{\operatorname{dist}}

\newcommand{\intl}{\int\limits}

\def\Xint#1{\mathchoice
    {\XXint\displaystyle\textstyle{#1}}%
    {\XXint\textstyle\scriptstyle{#1}}%
    {\XXint\scriptstyle\scriptscriptstyle{#1}}%
    {\XXint\scriptscriptstyle\scriptscriptstyle{#1}}%
    \!\int}
\def\XXint#1#2#3{\setbox0=\hbox{$#1{#2#3}{\int}$}
    \vcenter{\hbox{$#2#3$}}\kern-0.5\wd0}
\def\bint{\Xint-}
\def\dashint{\Xint{\raise4pt\hbox to7pt{\hrulefill}}}
\def\dashiint{\bint\kern-0.15cm\bint}
\newcommand{\ovl}[3]{\int_{#1}^{#2}\kern-#3pt\raise4pt\hbox to7pt{\hrulefill}\ }

\newcommand{\ovll}[3]{\intl_{#1}^{#2}\kern-#3pt\raise4pt\hbox to7pt{\hrulefill}\ }

\newcommand{\tvl}[2]{\iint_{#1}\kern-#2pt\raise4pt\hbox to7pt{\hrulefill}\ }

\newcommand{\bye}{\end{document}}

\newcommand{\tvls}[2]{\iint_{#1}\kern-#2pt\raise4pt\hbox to15pt{\hrulefill}\ }

\def\Q{\mathcal{Q}}
\def\K{\mathcal{K}}

 \def\B{\mathcal{B}}

\def\dist{\mathrm{dist}}

\def\P{\mathcal{P}}

\title[Liouville rigidity and time-extrinsic Harnack estimates] {Liouville rigidity and time-extrinsic Harnack estimates for an anisotropic slow diffusion}
\author[S. Ciani, U. Guarnotta
]{Simone Ciani \& Umberto Guarnotta
}
\address{Simone Ciani - Technische Universit\"at Darmstadt, Department of Mathematics,
Schlossgartenstr. 7,
64289, Darmstadt,
Germany} 
\address{Umberto Guarnotta - Universit\`a degli Studi di Palermo, Dipartimento di Matematica e Informatica, Via Archirafi 34, 90123, Palermo, Italy \vskip0.2cm}

\email{simone.ciani@unifi.it \& umberto.guarnotta@unipa.it}
\begin{document}
%%%%%%%%%%%%%%%%%%%%%%%%%%%%%%%%%%%%%%%%%%%%%%%%%%%
%%%%%%%%%%%%%%%%%%%%%%%%%%%%%%%%%%%%%%%%%%%%
\begin{abstract} 
We prove that non-negative solutions to the fully anisotropic equation 
\begin{equation*} 
    \partial_t u= \sum_{i=1}^N \partial_i (|\partial_i u|^{p_i-2} \partial_i u), \quad \qquad \text{in} \, \,  \R^N\times (-\infty, T), 
    % \quad \text{in} \quad \Omega_T= \Omega \times (0,T], \quad \text{with} \quad \Omega \subset \subset \R^N, 
\end{equation*} \noindent are constant if they satisfy a condition of finite speed of propagation and if they are both one-sided bounded, and bounded in $\R^N$ at a single time level. A similar statement is valid when the bound is given at a single space point. As a general paradigm, local H\"older estimates provide the basics for rigidity. Finally, we show that recent intrinsic Harnack estimates can be improved to a Harnack inequality valid for non-intrinsic times. Locally, they are equivalent.   


%%%%%%%%%%%%%%%%%%%%%%%%%%%%%%%%%%%%%%%%%%%%%%%%%%%

%%%%%%%%%%%%%%%%%%%%%%%%%%%%%%%%%%%%%%%%%%%%%%%%%%%
\end{abstract}

\maketitle

	\begin{center}
		\begin{minipage}{9cm}
			\small
			\tableofcontents
		\end{minipage}
	\end{center}

%%%%%%%%%%%%%%%%%%%%%%%%%%%%%%%%%%%%%%%%%%%%%%%%%%%%%%%%%%%%%%%%%%%%%%%
\section{Introduction to the problem}
\noindent 
Consider $u(x,t)$ as a function describing the temperature at the time $t$ of a point $x$ in an infinite isolated rod, being hence a solution of the heat equation. As usual, it is assumed that heat has spread from hotter zones to colder ones. Now, if one considers a non-negative solution in $\R^N \times (-\infty,0)$, the diffusive process has already gone on for an infinite amount of time, and it is reasonable to question if $u(x,t)$ has become constant. This fact, stated in this way, is generally false, as shown by the following examples:
\begin{equation}\label{examples}
    u_1(x,t)= e^{x_N+t}, \quad u_2(x,t)=e^{-t} \sin(x_1), \quad x \in \R^N.
\end{equation}
The two functions above are {\it eternal} solutions of the heat equation, i.e. solutions in $\R^N \times \R$. We call {\it ancient} solutions those solutions that solve the parabolic equation in $ \R^N \times (-\infty, T)$ for some time $T \in \R$. In line with the literature, we call {\it Liouville property} any rigidity condition that ensures the triviality of solutions. It is clear from $u_1$ that a sign condition is not enough to confirm our suspect, while the sign-changing solution $u_2$ shows that boundedness at a fixed time is not enough. Although Appel \cite{Appel} already proved in 1892 that an ancient solution to the heat equation which is two-sided bounded (as for instance $0 \leq u(x,t) \leq M$) is constant, the first optimal parabolic Liouville theorem for ancient solutions was found in 1952 by Hirschman (see \cite{Hirschman}, Bear \cite{Bear} and Widder \cite{Widder}, \cite{WidderBook} for the case $N=1$), stating that a non-negative ancient solution to the heat equation is constant if one adds the assumption that, for a time $t_o<T$,
\begin{equation}\label{log}
    \lim_{r \uparrow + \infty} \frac{\log(\sup_{|x|<r}u(x,t_o))}{r} \leq 0, \qquad \text{that is,} \qquad u(x,t_o) \leq e^{o(|x|)} \quad \text{as} \quad |x|\rightarrow +\infty.
\end{equation} This result was sharp in the sense that any function of the kind
\begin{equation*}\label{count-Hirsch}
    u(x,t)= e^{a^2t} \cosh{ax}
\end{equation*} shows that if \eqref{log} above is violated then $u(x,t)$ is not necessarily constant; but \eqref{log} is just a condition on the space variables for a fixed time. Sub-exponential optimal growth conditions have been generalized to different metric contexts; see for instance \cite{Sunra} and references therein. Not much later, in 1958, Friedman gave a condition on the behavior of non-negative ancient solutions to more general second-order parabolic equations as
\begin{equation}\label{seond-order}
    \partial_t u(x,t) = \sum_{i,j=1}^N a_{i,j}\frac{\partial^2 u(x,t)}{\partial x_i \partial x_j} + \sum_{i=1}^N b_i \frac{\partial_i u(x,t)}{x_i} +c u(x,t),  \quad \text{in} \quad \R^N \times \R_+,
\end{equation} being $b_i,c$ real numbers and $\{a_{i,j}\}_{i,j}$ a positive matrix. Now the assumption concerns infinite past times as
\begin{equation} \label{fried}
\lim_{t\downarrow - \infty} \frac{\log(u(0,t))}{t}= c+ \gamma, \qquad \gamma>0,\quad  c+\gamma \ge 0. \end{equation} \noindent See \cite{Friedman} for the result and \cite{Eidelman} for the earlier case of systems. Furthermore, conditions guaranteeing the stabilization of the solution to a constant were studied for a fixed space variable (see \cite{Eidelman-Kamin-Tedeev} and its references for an account). This short preamble is just to highlight that different assumptions, mainly on the second bound, may be requested to solutions of these parabolic equations in order to ensure Liouville property; it is therefore an incomplete list. The literature on these rigidity results is wide, so we refer the reader to the book \cite{QS-libro} and the survey \cite{Kogoj} for a more complete account.\vskip0.1cm
\noindent 
The heat equation can be regarded as a special case of the anisotropic $p$-Laplacian equation

\noindent \begin{equation} \label{EQ}
    \partial_t u= \sum_{i=1}^N \partial_i (|\partial_i u|^{p_i-2} \partial_i u),  
    % \quad \text{in} \quad \Omega_T= \Omega \times (0,T], \quad \text{with} \quad \Omega \subset \subset \R^N, 
\end{equation} \noindent
when $p_i\equiv 2$ for all $i=1,\dots,N$. When $2<p_i< \bar{p}(1+1/N)$, this equation describes the effect of competing diffusions along coordinate axes in finite speed of propagation (see \cite{Ant-Sh}, \cite{Mosconi}, \cite{Ruzicka} for an introduction to the parabolic problem). In one spatial dimension, the equation \eqref{EQ} is a one-dimensional $p$-Laplacian, which has a very interesting change of behavior from the degenerate case ($p>2$) to the singular one ($1<p<2$). Roughly speaking, the solutions to $p$-Laplacian singular equations behave more like solutions to elliptic equations and only one bound is enough to infer that they are constant (see \cite{DGV-Liouville} for more details). On the other hand, in the degenerate case two bounds are required to infer a Liouville property; see Section \ref{LiouvilleSection} below for a counterexample. In this paper we show some Liouville properties for non-negative solutions to \eqref{EQ} for a range of $p_i$s which is degenerate and allows a finite speed of propagation. In many physical circumstances, this is a more reasonable assumption than the sudden infinite expansion of the support of solutions to the heat equation.\vskip0.2cm 
\noindent The theory of regularity for solutions to \eqref{EQ}, even if much investigated, is still incomplete and fragmented (see, e.g., \cite{MingRadu} and \cite{Mar-survey} for an account on the elliptic case). The Liouville properties that we are about to describe are entailed by recent Harnack estimates, obtained with an approach of expansion of positivity. This has been shown relying on the behavior of abstract fundamental solutions in \cite{Ciani-Mosconi-Vespri}. Here we start from the aforementioned Harnack inequality (see Section \ref{Preliminaries}), which is formulated in an intrinsic geometry (see section Notations below) reflecting the natural scaling of the equation, and study some rigidity connections between local and global behavior of solutions.\newline
Similarly to the Liouville property inferred by Hirschman, we will prove that it is sufficient to have a one-sided bound (say, from below) and an estimate from the other side (say, from above), just for a fixed time. If these conditions are met, solutions are forced to be constant (Theorem \ref{Liouville1}). This clearly implies that a solution that is bounded both from above and below is constant; on the other hand, it is unreasonable to expect that just a one-sided bound suffices (see the example in Section \ref{LiouvilleSection}) for our range of $p_i$s. As a known fact, we comment that a precise decay on oscillation given by H\"older continuity estimates is enough as Liouville property, see Theorem \ref{Cutilisci}. This decay is usually easier to show than a complete Harnack inequality, and as such deserves its own attention: already in the range $\bar{p}(1+1/N)\leq p_i<\bar{p}(1+2/N)$, although continuity is expected by the regularizing properties of diffusion, no Harnack estimate may be available, because the competition among the diffusions is too strong (see for instance \cite{AS}). It would be of interest to compare this behavior to the one in porous materials by the sole control on the oscillation, for example as in \cite{Eurica}.\newline
On a similar track to Friedman's result, we prove that, fixed any spatial point and assuming that the solution is bounded at infinity in time, again it is forced to be constant (Theorem \ref{Liouville2}). Finally we state a Harnack inequality that frees the time variable to be intrinsic (Theorem \ref{WeakHarnack}, see \cite{DB} for the isotropic counterpart): being the Harnack estimate not anymore intrinsic in time, the estimate is more suitable for an application to rigidity. Moreover, this turns out to be useful to determine the optimal growth on the initial data when $|x|\rightarrow \infty $ for the solvability of the Cauchy problem for \eqref{EQ} (see for instance \cite{DB-He}). Clearly, this implies that the domain where the equation is solved must be, in turn, %intrinsic and
`compatible' with the anisotropy of the diffusion: this is certainly the case for ancient solutions.

\subsection*{Structure of the paper}
Section \ref{Preliminaries} is devoted to set up the functional framework and to recall some known properties of the solutions, as the existence of fundamental solutions, comparison principles, and the Harnack inequality. 
Section \ref{Appendix} is concerned with the study of H\"older continuity of solutions. In Section \ref{LiouvilleSection} we prove the Liouville-type results, while Section \ref{WHSection} pertains an alternative formulation of the Harnack inequality, which turns out to be locally equivalent to the known one.

\section*{Notations} 
\begin{itemize}
\item[-] We define the following function of $p_i$s, called {\it harmonic mean}: $\,
\bar{p}=N(\sum_{i=1}^N1/p_i)^{-1}$.\newline We suppose that $ 1<p_1\leq p_2 \leq \ldots \leq p_N$ are ordered, as well as $\bar{p}<N$.
\vskip0.1cm \noindent 
    
\item[-] For any $\rho, \theta>0$ and $x\in\R^N$, we denote by $K_{\rho}(x) \subset \R^N$ the cube of side $2\rho$ centered at $x$.
% Analogously, $B_{\rho}^d(x)\subset \R^d$ denotes the $d$-dimensional open ball of radius $\rho$ and center $x$. 
\noindent Let $x_o+\K_{\rho}(\theta)$ stand for the anisotropic cube of radius $\rho$, ``magnitude'' $\theta$, and center $x_o$, i.e.,
\begin{equation}\label{anisocubi}
x_o+\K_{\rho}(\theta)= \prod_{i=1}^N\bigg{\{}|x-x_{o,i}|<\theta^{{(p_i-\bar{p})}/{p_i}}\rho^{{\bar{p}}/{p_i}}\bigg{\}}.
\end{equation}
If either $\theta=\rho$ or $p_i=p$ for all $i=1,\ldots,N$, then $x_o+\K_{\rho}(\theta)=K_{\rho}(x_o)$.

\item[-] For any $\rho, \theta,C >0$ and $(x_o,t_o) \in \R^{N+1}$, we consider the following anisotropic cylinders:
\begin{equation*}\label{cylinders}
\begin{cases}
\text{centered: }(x_o,t_o)+\Q_{\rho}(\theta,C)=
(x_o+\K_{\rho}(\theta) )\times (t_o-\theta^{2-\bar{p}}(C\rho)^{\bar{p}},t_o+\theta^{2-\bar{p}}(C\rho)^{\bar{p}});\\
\text{forward: }(x_o,t_o)+\Q^+_{\rho}(\theta,C)= (x_o+\K_{\rho}(\theta) )\times [t_o,t_o+\theta^{2-\bar{p}}(C\rho)^{\bar{p}});\\
\text{backward: }(x_o,t_o)+\Q^-_{\rho}(\theta,C)=
(x_o+\K_{\rho}(\theta) )\times (t_o-\theta^{2-\bar{p}}(C\rho)^{\bar{p}},t_o].
\end{cases}
\end{equation*} \noindent We omit the index $C$ when the constant is clear from the context. \vskip0.1cm \noindent 
\item[-] For $\Omega \subset \subset \R^N$, i.e., $\Omega$ open and bounded set in $\R^N$, we denote with $\Omega_T= \Omega \times [-T,T]$, $T>0$, the parabolic domain, and with $S_s= \R^N \times (-\infty, s)$, $s\in \R$, the space strip.\vskip0.1cm \noindent 
\item[-] We adopt the convention that the constant $\gamma>0$ may change from line to line, when depending only on fixed quantities $\{N,p_i\}$.
\end{itemize}

\section{Preliminaries and Tools of the Trade} \label{Preliminaries}
% Consider the following equation in $\Omega_T=[0,T] \times \Omega$, $\Omega \subset \subset \R^N$ open,
% \begin{equation} \label{prototype}
%     u_t= \sum_{i=1}^N (|u_{x_i}|^{p_i-2} u_{x_i})_{x_i}, \quad p_i>2 \, \quad \forall i=1,..N.
% \end{equation} \noindent 
\noindent 
We begin with the definition of solution. For $\Omega \subseteq \R^N$ open rectangular domain and $T>0$, we set $\Omega_T= \Omega \times [-T,T]$ and define the Banach spaces 
% \[ W^{1,{\bf{p}}}_o(\Omega):= \{ u \in W^{1,1}_o(\Omega) |\,  \partial_i u \in L^{p_i}(\Omega) \}, \]
\[W^{1,{\bf{p}}}_{loc}(\Omega):= \{ u \in W^{1,1}_{loc}(\Omega) |\,  \partial_i u \in L^{p_i}_{loc}(\Omega) \}, \]
% \[ L^{{\bf{p}}}(0,T;W^{1,{\bf{p}}}_o(\Omega)):= \{u \in L^1(0,T;W^{1,1}_o(\Omega))|\, \partial_i u \in L^{p_i}(0,T;L^{p_i}_{loc}(\Omega))   \}, \]
\[ L^{{\bf{p}}}_{loc}(0,T;W^{1,{\bf{p}}}_{loc}(\Omega)):= \{u \in W^{1,1}_{loc}(0,T;L^1_{loc}(\Omega))|\, \partial_i u \in L^{p_i}_{loc}(0,T;L^{p_i}_{loc}(\Omega))   \}. \]
%We use the symbol $u_{x_i}$ meaning the $i$-th weak derivative for a function $u \in W^{1,\bf{p}}_0(\Omega)$, and the partial derivatives $\frac{\partial u}{\partial x_i}$ for formal calculations.\newline
These are usually called anisotropic spaces (see for instance \cite{Ant-Sh}). When $\bar{p}>N$ and $\partial \Omega$ is regular enough, the space $W^{1,{\bf{p}}}(\Omega)$ is embedded in the space of H\"older continuous functions \cite{VenTuan}. A function \[ u \in C_{loc}(0,T; L^2_{loc}(\R^N)) \cap L^{\bf{p}}_{loc}(0,T;W^{1,{\bf{p}}}_{loc}(\R^N))\] is called a {\it local weak solution} of \eqref{EQ} in $S_T$ if, for all $0<t_1<t_2<T$ and any $\varphi \in C^{\infty}_{loc}(0,T;C_o^{\infty}(\R^N))$,
\begin{equation} \label{anisotropic-localweaksolution}
\int_{\R^N} u \varphi \, dx \bigg|_{t_1}^{t_2}+ \int_{t_1}^{t_2} \int_{\R^N} (-u \, \varphi_t + \sum_{i=1}^N \, |\partial_i u|^{p_i-2} \partial_i u \, \partial_i \varphi) \, dx dt=0,
\end{equation} \noindent Similarly, when considering $\Omega\subset \R^N$ bounded set, by a {\it local weak solution} to \eqref{EQ} in $\Omega_T$ we mean a function $u \in C_{loc}(0,T; L^2_{loc}(\Omega)) \cap L^{\bf{p}}_{loc}(0,T;W^{1,{\bf{p}}}_{loc}(\Omega))$ satisfying for all compact sets $K \subset \Omega$ and for all $\varphi \in C^{\infty}_{loc}(0,T;C_o^{\infty}(K))$ the integral equality
\begin{equation} \label{localweaksolution}
\int_{K} u \varphi \, dx \bigg|_{t_1}^{t_2}+ \int_{t_1}^{t_2} \int_{K} (-u \, \varphi_t + \sum_{i=1}^N \, |\partial_i u|^{p_i-2} \partial_i u \, \partial_i \varphi) \, dx dt=0,\quad \text{for all} \quad 0<t_1<t_2<T.
\end{equation}
% By a density and approximation argument this actually holds for any test function $\varphi$ belonging to \[ W^{1,2}_{loc}(0,T;L^2_{loc}(\R^N))\cap L^{\bf{p}}_{loc}(0,T;W^{1,{\bf{p}}}_o(\Omega)),\] provided $ \Omega \subset \subset \R^N $ is a rectangular domain (see Theorem 3 in \cite{Haskovec-Schmeiser}). Now, by a local weak solution\footnote{\textcolor{red}{Sembra che la definizione non sia univoca; vedi sopra.}} to the equation \eqref{EQ} in $\R^N\times \R_+$, we mean a function
% \[ u \in L^2_{loc}(\R; L^2_{loc}(\R^N)) \cap L^{\bf{p}}_{loc}(\R;W^{1,{\bf{p}}}_{loc}(\R^N))\] that, for each choice of $T\in \R$ and $\Omega \subset \subset \R^N$, is a weak local solution in $\Omega_T$ in the sense of the previous definition. The obvious adjustments define the concept of local weak solution to the equation in $S_T$, $T\in \R$.
Now we briefly introduce the main tools for our proofs: the intrinsic Harnack inequality, the existence of an abstract Barenblatt-type solution, and a local comparison principle.\newline
Hereafter, with the only exception of Theorem \ref{Cutilisci}, we restrict our attention to the range
\begin{equation}\label{pi}
    2<p_1\leq p_N<\bar{p}(1+1/N),\qquad \bar{p}<N,
\end{equation}
and we will refer to the constants $C_i$, $i=1,2,3$, appearing in the following theorem.
\begin{theorem} \label{Harnack-Inequality}
Let $u\ge 0$ be a  local weak solution to \eqref{EQ} in $\Omega_T$ and let \eqref{pi} be valid. Suppose that $u(x_o,t_o)>0$ for a Lebesgue point $(x_o,t_o) \in \Omega_T$ for $u$. Then there exist $C_{1}\ge 0, C_3\ge C_2\ge 1$, depending only on $N$ and the $p_{i}$s, such that, letting $\theta=u(x_o,t_o)/C_1$, it holds
\begin{equation}\label{Harnack}
  \frac{1}{C_{3}}\sup_{x_o+\K_{\rho}(\theta)}u(\,\cdot\, , t_o - \theta^{2-\bar p}\, (C_{2}\, \rho)^{\bar p} )\le  u(x_o,t_o) \le C_{3} \inf_{x_o+\K_{\rho}(\theta)} u(\,\cdot\, ,   t_o + \theta^{2-\bar{p}}\, (C_{2}\, \rho)^{\bar{p}})
    \end{equation}
    with $\K_{\rho}(\theta)$ defined as in \eqref{anisocubi}, whenever $\rho, \theta>0$ satisfy 
    \begin{equation} \label{side-condition}
    \theta^{2-\bar p}\, (C_{3}\,  \rho)^{\bar p}<T-|t_o| \qquad \text{and} \qquad x_o+\K_{C_{3}\, \rho}(\theta)\subseteq \Omega.\end{equation}
\end{theorem}

\noindent The assumption $u(x_o,t_o)>0$ is understood by a suitable limit process, as customary. Semi-continuity clarifies this definition, as long as a theoretical maximum principle is in force (see \cite{CianiGuarnotta}, \cite{Mosconi}, \cite{Liao} for an account). 
% More precisely, $u(x_o,t_o)>0$ is interpreted on an integral average on classical balls, i.e. 
% \begin{equation}\label{pointsense}
% \lim_{r\downarrow 0} \bigg( \dashint_{B^{N+1}_r(x_o,t_o)} u(x,t)\, dxdt\bigg)>0.
% \end{equation} \noindent Theorem \ref{Harnack-Inequality} has an easier formulation if solutions are supposed to be continuous. 
Theorem \ref{Harnack-Inequality} has been proved in \cite{Ciani-Mosconi-Vespri} without the assumption of H\"older continuity of solutions, which can be shown (see Section \ref{Appendix}) to be a sole consequence of \eqref{Harnack}. This important property has been faced several times in the past, with imprecise proofs or an unclear geometric setting. For this reason, and in order to explain the main adversities that anisotropic diffusion obliges us to face, we include in Section \ref{Appendix} a proof of local H\"older continuity of solutions to \eqref{EQ}, which follows Moser's ideas \cite{Moser} through an appropriate anisotropic intrinsic geometry. Taking for granted their continuity, in what follows we will refer directly to the point-wise values of solutions.\vskip0.1cm \noindent Let us comment Theorem \ref{Harnack-Inequality} from a global point of view: if we pick a point $(x_o,t_o) \in \Omega_T$ where $u$ is positive, it is possible to `detect' the sets where the pointwise controls \eqref{Harnack} hold true. This is the core of the next proposition.

\begin{proposition}\label{paraboloids}
Suppose the assumptions of Theorem \ref{Harnack-Inequality} are satisfied for $(x_o,t_o) \in \Omega_T$. Then
\begin{equation} \label{estimate-paraboloid}
\inf_{\P_\theta^+(x_o,t_o)} u \ge u(x_o,t_o)/C_3 \qquad \text{and}\qquad \sup_{\P_\theta^-(x_o,t_o)}u \leq C_3 u(x_o,t_o),\end{equation}
where, setting $\theta= u(x_o,t_o)/C_1$, the paraboloids $\P^+_{\theta}(x_o,t_o)$ and $\P^-_{\theta}(x_o,t_o)$ are defined by
\[
\P_\theta^+(x_o,t_o)= \bigg{\{}(x,t) \in \Omega_T:\, \,  C_2^{\bar{p}} |x_i-x_{o,i}|^{p_i}\theta^{2-p_i}\leq  (t-t_o)\leq C_2^{\bar{p}}\varrho^{\bar{p}}\theta^{2-\bar{p}}, \, \,   \forall i=1,..N \bigg{\}},
\]
\[
\P_\theta^-(x_o,t_o)= \bigg{\{}(x,t) \in \Omega_T:\,\,    -C_2^{\bar{p}}\varrho^{\bar{p}}\theta^{2-\bar{p}}\leq (t-t_o) \leq -C_2^{\bar{p}}|x_i-x_{o,i}|^{p_i} {\theta}^{2-p_i},  \, \,   \forall i=1,..N \bigg{\}},
\] 
% The union of these two sets through the point $(x_o,t_o)$ is the whole anisotropic intrinsic paraboloid
% \begin{equation} \label{paraboloide}
%     \P_u(x_0,t_0)= \bigg{\{}(x,t) \in \Omega_T:\quad  |x_i-x_{0,i}|^{p_i}\leq C_2^{-\bar{p}} \theta^{p_i-2} |t-t_0|\leq (\rho^{+})^{\bar{p}}\theta^{p_i-\bar{p}}, \quad \forall i=1,..N \bigg{\}},
%     % = \P_{\theta}^+(x_0,t_0) \cup \P_{\theta}^-(x_0,t_0) .
% \end{equation}
\noindent with $\varrho$ depending on $u$, $\Omega_T$, and $(x_o,t_o)$ according to the following expression: 
\begin{equation}\label{rho+}
    \varrho^{\bar{p}}= C_3^{-\bar{p}} \bigg( \frac{u(x_o,t_o)}{C_1} \bigg)^{\bar{p}-2} \min_{i=1,\dots,N} \bigg{\{}(T-|t_o|), \, \bigg( \frac{\dist(x_o, \partial \Omega)}{2}\bigg)^{p_i} \bigg( \frac{u(x_o,t_o)}{C_1} \bigg)^{2-p_i} \bigg{\}}.
\end{equation}
\end{proposition}
\noindent It is remarkable that estimate \eqref{Harnack} is prescribed on a {\it{space}} configuration depending on the solution, in contrast to what happens with $p$-Laplacian type equations. This is due to the natural scaling of the equation (see \cite{CianiGuarnotta}), because the expansion of positivity of solutions is readily checked via comparison with the following family of Barenblatt-type solutions.

\begin{theorem}
\label{Barenblatt}
Set $\lambda=N(\bar{p}-2)+\bar{p}$ and let \eqref{pi} be satisfied. For each $\sigma >0$ there exists $\tilde{\eta}>0$ and a local weak solution $\B_{\sigma}(x, t)$ to \eqref{EQ} 
with the following properties, valid for any $t\in(0,T)$:
\begin{enumerate}
\item $\displaystyle{\|\B_{\sigma}(\cdot, t)\|_{\infty}=\sigma \,  t^{-\alpha}}$,
\item
$\displaystyle{{\rm supp}(\B_{\sigma}(\cdot, t))\subseteq \prod_{i=1}^N \big{\{} |x_i|\le    \sigma^{(p_i-2)/p_i}\, t^{\alpha_i} \big{\}}}$, $\qquad \qquad$  $\alpha=N/\lambda$, $\alpha_i=(1+2\alpha)/p_i-\alpha$, 
\item
$\displaystyle{\{\B_{\sigma}(\cdot, t)\ge \eta\,  \sigma \, t^{-\alpha}\}\supseteq \prod_{i=1}^N \big{\{} |x_i|\le \eta\, \sigma^{(p_{i}-2)/p_{i}}\, t^{\alpha_i} \big{\}}=:\mathcal{P}_t}$.
\end{enumerate}
\end{theorem}
\noindent The existence of a Barenblatt Fundamental solution $\B$ is a consequence of the finite speed of propagation of solutions to \eqref{EQ} combined with a particular correspondence of the Cauchy problems associated to \eqref{EQ} and to an anisotropic Fokker-Planck equation. On the other hand, the properties of $\B$ stated above stem from comparison techniques and the invariance of the equation \eqref{EQ} under scaling, which entitles $\B$ to be a self-similar solution. We refer to \cite{Ciani-Mosconi-Vespri} for the proofs of these facts and the following proposition; see also \cite{CSV}, \cite{Vazquez} for the singular case. 
\begin{proposition} \label{local-comparison}
Let $\Omega \subset \R^N$ be a bounded open set and $u,v$ be local weak solutions to the equation \eqref{EQ} in $\Omega_T$. Let $\tilde{\Omega} \subset \Omega$ and $0<\tilde{T}<T$. If $u,v$ satisfy $u(x,t) \ge v(x,t)$ in the parabolic boundary of $\tilde{\Omega}_{\tilde{T}}$, then $u \ge v$ in $\tilde{\Omega}_{\tilde{T}}$. 
\end{proposition} \noindent The point-wise boundary inequality assumed in Proposition \ref{local-comparison} will be used in the proof of Theorem \ref{WeakHarnack} locally, and as such, it has a well-defined meaning thanks to the results of the next section.

\section{H\"older Continuity of solutions}
\label{Appendix}

\begin{theorem}\label{HC}
Under condition \eqref{pi}, any local weak solution $u$ to \eqref{EQ} is locally H\"older continuous. More precisely, there exist $\gamma>1$ and $\chi \in (0,1)$, depending only upon $p_i,N$, with the following property: for each compact set $K \subset \subset \Omega_T$ there exist a set $\Lambda$ and $\omega_o=\omega_o(K, \|u\|_{\infty,K})$ such that $K \subset \Lambda \subseteq \Omega_T$ and, for every $(x,t)$, $(y,s) \in K$,
\begin{equation} \label{HContinuity}
    |u(x,t)-u(y,s)| \leq \gamma \omega_o \bigg(\frac{\sum_{i=1}^N |x_i-y_i|^{{p_i}/{\bar{p}}}\omega_o^{{(\bar{p}-p_i)}/{\bar{p}}}+ |t-s|^{1/{\bar{p}}}\omega_o^{{(\bar{p}-2)}/{\bar{p}}}}{{\bf{p}}\text{-dist}(K,\partial \Lambda) } \bigg)^{\chi},
\end{equation}\noindent with
\begin{equation} \label{pi-dist}
\begin{aligned}
&{\bf{p}}\text{-dist}(K,\partial \Lambda):=\inf \{ {\bf{p}}_x, {\bf{p}}_t \}, \quad \text{being}\\
& {\bf{p}}_x=\inf \bigg{\{} 
|x_i-y_i|^{{p_i}/{\bar{p}}}(\omega_o/C_1)^{{(\bar{p}-p_i)}/{\bar{p}}}\, : \, (x,t) \in K, (y,s) \in \partial \Lambda,\, i=1,..,N\bigg{\}},\\
& {\bf{p}}_t=\inf \bigg{\{} 
|t-s|^{{1}/{\bar{p}}}(\omega_o/C_1)^{{(\bar{p}-2)}/{\bar{p}}}\, : \, (x,t) \in K, (y,s) \in \partial \Lambda\bigg{\}}.
\end{aligned}
\end{equation} \noindent Furthermore, if $u$ is bounded in $\Omega_T$ then \eqref{HContinuity} holds with $\Lambda= \Omega_T$.
\end{theorem} \noindent We prove Theorem \ref{HC} in four steps, without assuming that $u$ is globally bounded.
\begin{proof}
Let us fix a compact set $K \subset \subset \Omega_T$ and two points $(y,s), (x,t) \in K$. \vskip0.2cm 
\noindent {\small{STEP 1-{\it A global bound for the solution in $K$.}}}
\vskip0.2cm \noindent Let $\bar{p}_2= \bar{p}(1+2/N)$ and for $k>0$ we define the increasing functions $g(k)=\sum_{i=1}^{N}k^{ p_{i}-2}$ and $h(k)=\left(\sum_{i=1}^{N}k^{p_{i}-\bar p_{2}}\right)^{-1}$. We use the estimates in \cite[Lemma 4.2]{Mosconi}: under condition \eqref{pi}, there exists $\tilde{\gamma}>0$ such that solutions to \eqref{EQ} satisfy
\begin{equation}
\label{supest}
\|u_{+}\|_{L^{\infty}(Q_{\lambda/2, M})}\leq g^{-1}(1/M)+ h^{-1}\left({\gamma}\Big(M\, \dashiint_{Q_{\lambda, M}} u_{+}^{\bar p_{2}}\, dx\Big)^{{\bar p}/{(N+\bar p)}}\right), 
\end{equation} in the (non-intrinsic) anisotropic cylinders 
\begin{equation} \label{anisocylinder}
Q_{\lambda, M}= \prod_{i=1}^{N}\left[-\lambda^{{1}/{p_{i}}}, \lambda^{{1}/{p_{i}}}\right]\times [-M\, \lambda, 0],\quad \quad  M, \lambda>0.
\end{equation}
\vskip0.2cm \noindent By compactness of $K$, we find $(x_i,t_i) \in K$ and $\lambda_i, M_i \in \mathbb{R}_+$, $i=1,\dots,m$, for $m\in\N$, such that 
\begin{equation*}
    K \subset \Lambda:=\bigcup_{j=1}^m \{(x_j,t_j)+Q_{\lambda_j,M_j}\}\subseteq \bigcup_{j=1}^m \{(x_j,t_j)+Q_{2\lambda_j,M_j} \}\subseteq \Omega_T,
\end{equation*} \noindent being $Q_{\lambda,M}$ as in \eqref{anisocylinder}. 
% Indeed, $D=\dist(K, \partial \Omega_T)>0$, so we may consider $\lambda_j,M_j< \min\{(D/2)^{p_i},\sqrt{D/2} \}$.
\noindent According to \eqref{supest}, for each anisotropic cylinder  $\hat{Q}_{\lambda_j,M_j}=(x_j,t_j)+Q_{\lambda_j,M_j}$, $j=1,\dots,m$, we deduce the estimate 
\begin{equation*} \begin{aligned}\label{A}
\| u\|_{L^{\infty}(\hat{Q}_{\lambda_j, M_j})} &
% \leq g^{-1} (1/M_j)+ h^{-1} \bigg( C \bigg(M_j \dashint \dashint_{\hat{Q}_{2\lambda_j,M_j}} |u|^{\bar{p_2}}\,  dxdt\bigg)^{\frac{\bar{p}}{N+\bar{p}}} \bigg)\\
% &
\leq g^{-1} (1/\min_{j} M_j)+ h^{-1} \bigg( \gamma \max_{j=1,\dots,m}  \bigg( M_j \dashint \dashint_{\hat{Q}_{2\lambda_j,M_j}} |u|^{\bar{p}_2}\,  dxdt\bigg)^{{\bar{p}}/{(N+\bar{p})}} \bigg)=: \mathcal{I},
\end{aligned}
\end{equation*} \noindent because $h$,$g$, are monotone increasing.
% Observe that $\mathcal{I}$ does not depend on $j\in \{1,\dots,m\}$.
\noindent Finally, we define $\omega_o=\omega_o(K)$ as \begin{equation} \label{0}
\omega_o:= 2\mathcal{I}. \end{equation} \noindent  Accordingly,
\[
K \subset \bigcup_{j=1}^m \hat{Q}_{\lambda_j, M_j}(x_j,t_j)= \Lambda \qquad \mbox{and} \qquad 2 \|u\|_{L^{\infty}(\Lambda)} \leq \omega_o.
\]
\vskip0.2cm \noindent {\small{STEP 2-{\it Accommodation of degeneracy and alternatives.}}}
\vskip0.2cm \noindent
Recalling \eqref{pi-dist}
% \begin{equation} \label{pi-dist}
% \begin{aligned}
% &{\bf{p}}\text{-dist}(K,\partial \Omega_T):=\inf \{ {\bf{p}}_x, {\bf{p}}_t \}, \quad \text{being}\\
% & {\bf{p}}_x=\inf \bigg{\{} 
% |x_i-y_i|^{\frac{p_i}{\bar{p}}}(\omega_o(K)/C_1)^{\frac{\bar{p}-p_i}{\bar{p}}}\, : \, (x,t) \in K, (y,s) \in \partial \Omega_T,\, i=1,..,N\bigg{\}},\\
% & {\bf{p}}_t=\inf \bigg{\{} 
% |t-s|^{\frac{1}{\bar{p}}}(\omega_o(K)/C_1)^{\frac{\bar{p}-2}{\bar{p}}}\, : \, (x,t) \in K, (y,s) \in \partial \Omega_T\bigg{\}},
% \end{aligned}
% \end{equation} \noindent
we define $R:= [{\bf{p}}\text{-dist}(K,\partial \Lambda)]/(2C_3)$. Now, by definition of $R$, the intrinsic cylinder centered at $(y,s)\in K$ and constructed with $R$ and $\omega_o$ is contained inside  $\Lambda$, that is,
\[
(y,s)+ \Q_{R}(\omega_o/C_1,C_2) \subseteq \Lambda.
\]
\noindent Now consider any other point $(x,t) \in K$. We reduce the study of the oscillation only in $(y,s) + \Q_R^-(\omega_o/C_1,C_2)$, having elsewhere the H\"older continuity of $u$. Indeed, if $|s-t| \ge (\omega_o/C_1)^{2-\bar{p}}(C_2 R)^{\bar{p}}$, we have
\[
|u(y,s)-u(x,t)|\leq |u(y,s)|+|u(x,t)|\leq \omega_o \leq 2C_3 \omega_o \bigg( \frac{(\omega_o/C_1)^{{(\bar{p}-2)}/{\bar{p}}}|s-t|^{{1}/{\bar{p}}}}{{\bf{p}}\text{-dist}(K,\partial \Lambda)} \bigg) \] by definition of $R$.
%\[
%{\bf{p}}\text{-dist}(K,\partial \Lambda) \leq 2C_3 (\omega_o/C_1)^{{(\bar{p}-2)}/{\bar{p}}} |s-t|^{{1}/{\bar{p}}}
%\qquad \Rightarrow \qquad
%1\leq 2C_3 \bigg( \frac{(\omega_o/C_1)^{{(\bar{p}-2)}/{\bar{p}}}|s-t|^{{1}/{\bar{p}}}}{{\bf{p}}\text{-dist}(K,\partial \Lambda)} \bigg).
%\]
Similarly, if $|y_i-x_i| \ge (\omega_o/C_1)^{{(p_i-\bar{p})}/{\bar{p}}} R^{{\bar{p}}/{p_i}}$ for some $i \in \{1,\dots,N\}$, the same conclusion follows from
\[
|u(y,s)-u(x,t)|\leq |u(y,s)|+|u(x,t)|\leq \omega_o \leq 2C_3 \omega_o \bigg( \frac{(\omega_o/C_1)^{{(\bar{p}-p_i)}/{\bar{p}}}|y_i-x_i|^{{p_i}/{\bar{p}}}}{{\bf{p}}\text{-dist}(K,\partial \Lambda)} \bigg).
\]
This technical stratagem justifies the definition \eqref{pi-dist}. Hence we can assume that
\begin{equation}\label{exclusion}
   |s-t|<(\omega_0/C_1)^{2-\bar{p}} (C_2 R)^{\bar{p}} \quad \text{and}  \quad  |y_i-x_i|< (\omega_o/C_1)^{{(p_i-\bar{p})}/{p_i}} R^{{\bar{p}}/{p_i}}  \quad \forall i=1,\dots, N, 
\end{equation}
that is, 
\[(x,t) \in (y,s)+ \Q_R^-(\omega_o/C_1,C_2). \]
We take the cylinder $\Q_0:=(y,s)+\Q_R^-(\omega_o/C_1,C_2)$ as the first element of a net $\{\Q_n\}_n$ of cylinders shrinking to the center $(y,s)$. This net will be constructed to control uniformly the oscillation.

\vskip0.2cm \noindent

\noindent {\small{STEP 3-{\it Controlled reduction of oscillation.}}} 
% \vskip0.2cm \noindent
% The reduction of oscillation will be a consequence of the following fact
\begin{proposition}\label{birra} Let the hypothesis of Theorem \ref{HC} be valid and assume also \eqref{exclusion}. Then, setting 
\begin{equation*}
    \begin{cases}
    \omega_0= \omega_o(K),\\
    \omega_n=\delta \omega_{n-1}, \, n\ge 1,
    \end{cases} 
    \begin{cases}
    \theta_n= \omega_n/C_1, \, n\ge 0,\\
    \rho_0=R,\\
    \rho_n= \varepsilon \rho_{n-1}, \, n\ge 1,\\
    \end{cases}
    \begin{cases}
    \delta=4C_3/(1+4C_3),\\
    \varepsilon=\delta^{{(\bar{p}-2)}/{\bar{p}}}/A, \\
    A=4^{p_N}, \end{cases}
\end{equation*} \noindent we have both the inclusions
\[
\Q_{n}\subset \Q_{n-1}, \quad \text{with} \quad \Q_n= (y,s)+ \Q_{\rho_n}^-(\theta_n)= \prod_{i=1}^N \bigg{\{}|y_i-x_i|<\theta_n^{{(p_i-\bar{p})}/{p_i}}\rho_n^{{\bar{p}}/{p_i}}  \bigg{\}} \times \bigg(s-\theta_n^{2-\bar{p}} (C_2\rho_n)^{\bar{p}} ,\, s\bigg],
\] and the inequalities
\begin{equation}\label{control}
    \osc_{\Q_n} u \leq \omega_n = \delta^n \omega_o.
\end{equation}
\end{proposition}

\begin{proof}[Proof of Proposition \ref{birra}]\noindent First of all, we prove that $\Q_{n} \subset \Q_{n-1}$ for all $n\in\N$. By direct computation,
%The first step is $\osc_{\Q_{0}} u =\osc_{Q_R(\omega_o/C_1)} \leq \omega_o$, that holds true because $\omega_o \ge 2 \osc_{Q_0}$ having the inclusion $(y,s)+Q_R(\omega_o/C_1)\subset \Lambda$ by the previous accommodation of degeneracy. Moreover, we postpone the proof of $\Q_1\subseteq \Q_0=(y,s)+Q_R(\omega_o/C_1)$ for the following general fact.
\begin{equation*} \begin{aligned}
\theta_{n}^{2-\bar{p}} (C_2\rho_{n})^{\bar{p}} = \bigg(\frac{\delta \omega_{n-1}}{C_1}\bigg)^{2-\bar{p}}\bigg((C_2\rho_{n-1}/A)^{\bar{p}}\delta^{\bar{p}-2}\bigg)= \theta_{n-1}^{2-\bar{p}} (C_2\rho_{n-1}/A)^{\bar{p}}.
\end{aligned} \end{equation*}  For each $i\in \{1,..,N\}$, since $p_i>2$ and $\delta \in (0,1)$, it holds
\[
\theta_{n}^{p_i-\bar{p}} \rho_{n}^{{\bar{p}}} = \delta^{p_i-2} \theta_{n-1}^{{p_i-\bar{p}}} ( {\rho_{n-1}}/{A} )^{{\bar{p}}} \leq \theta_{n-1}^{{p_i-\bar{p}}} ( {\rho_{n-1}}/{A} )^{{\bar{p}}}.
\] This computation shows a little more, by allowing indeed $\Q_{n}\subset (y,s)+\Q_{\rho_{n-1}/A}^-(\theta_{n-1})\subset \Q_{n-1}$.
\noindent Now we prove \eqref{control} by induction. The base step holds true: indeed, the accommodation of degeneracy (see Step 2 above) entails $\Q_0\subset\Lambda$, so that the bound produced in Step 1 yields
\[
\osc_{\Q_0} u \leq \osc_{\Lambda} u \leq 2\|u\|_{L^\infty(\Lambda)} \leq \omega_o.
\]
We assume now that the statement \eqref{control} is true until step $n$ and we show it for $n+1$. This will determine the number $A$. More precisely, we assume that $\osc_{\Q_n}u \leq \omega_n$ and, by contradiction, that $\osc_{\Q_{n+1}} u > \omega_{n+1}$. We set
\[
M_{n}= \sup_{\Q_{n}} u, \qquad m_{n}= \inf_{\Q_{n}}u, \qquad P_{n}=(y,\, s-\theta_{n}^{2-\bar{p}}(C_2\rho_{n})^{\bar{p}}).
\] Now we observe that one of the following two inequalities must be valid:
\[
M_{n}-u(P_{n}) > \omega_{n+1}/4 \qquad \text{or} \qquad u(P_{n})-m_{n} > \omega_{n+1}/4.
\] Indeed, if both alternatives are violated, then by adding the opposite inequalities we obtain $\osc_{\Q_{n}} u  \leq \omega_{n+1}/2< \osc_{\Q_{n+1}}$, generating a contradiction with $\Q_{n+1} \subseteq \Q_n$. Let us suppose $M_{n}-u(P_{n}) \ge  \omega_{n+1}/4$, the other case being similar. In particular we have the double bound
\begin{equation}
    \label{doublebound}
    \omega_{n+1}/4\leq M_n - u(P_n) \leq \omega_n.
\end{equation}
Let us set $\hat{\theta}_n=(M_n-u(P_n))/C_1$. We work in the half-paraboloid $\P^+_n=\P^+_{\hat{\theta}_n}(P_n)$ for times restricted to the ones of $\Q_n$.
% , which is
% \begin{equation*}
% \begin{split}
% \mathcal{P}^+_n =\bigg{\{} |x_i-y_i|^{p_i}<  (C_2)^{-\bar{p}} [(M_n-u(P_n))/C_1]^{p_i-2}(t-s+& (C_2\rho_n)^{\bar{p}}(\omega_n/C_1)^{2-\bar{p}}), \\
% &s-(C_2\rho_n)^{\bar{p}}(\omega_n/C_1)^{2-\bar{p}})\leq t\leq s \bigg{\}}.
% \end{split}
% \end{equation*}
\noindent
% By \eqref{doublebound} we deduce $\hat{\theta}_n \leq \theta_n$, so the condition $\rho_n^{\bar{p}} \theta_n^{2-\bar{p}}<\varrho^{\bar{p}} \hat{\theta}_n^{2-\bar{p}}$ specified by the request $\P_n^+ \subset \Omega_T$ is satisfied (see Figure \ref{FigA}): indeed, recalling the definition of $R$ and $\hat{\theta}_n \leq \theta_n$, besides that $\P_n^+$ is restricted to the times of $\Q_{n}$, we have the following estimate along the time variable:
% \[\rho_n^{\bar{p}} \bigg(\frac{\omega_n}{C_1}\bigg)^{2-\bar{p}}\leq R^{\bar{p}} \bigg( \frac{\omega_o}{C_1}\bigg)^{2-\bar{p}}\leq  \varrho^{\bar{p}} \bigg( \frac{M_n-u(P_n)}{C_1}\bigg)^{2-\bar{p}}.\]
%  It is worth observing that, in the previous estimate, the number $\varrho$ depends on the function $M_n-u(P_n)$ and the point $P_n$.
The starting time of $P_n^+$ is the same as the one of $\Q_n$ (see Figure \ref{FigA}).

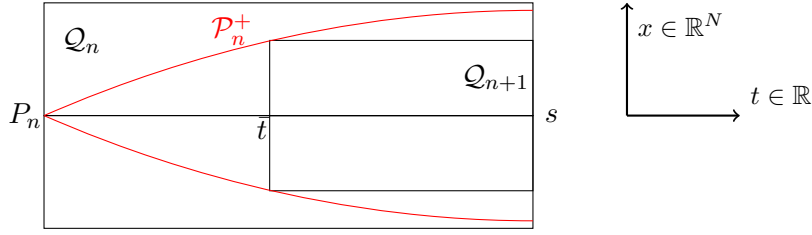
\begin{figure}
\centering
\begin{tikzpicture}[scale=0.25]

\draw[thick,->] (25,0) -- (25,6) node[anchor=north west] {\small{$x \in \R^N$}};
\draw[thick,->] (25,0) -- (31,0) node[anchor=south west] {\small{$t \in \R$}};
\draw (20,0) rectangle (-6,6);

\draw (20,0)  rectangle (6,4);

\draw (20,0) rectangle (-6,-6);

\draw (20,0) rectangle (6,-4);

\draw (18, 2) node{$\Q_{n+1}$};
\draw (-4, 4) node{$\Q_{n}$};
\draw (4, 4.6) node{\textcolor{red}{$\P^+_n$}};
\draw (-7, 0) node{$P_n$};
\draw (21, 0) node{$s$};

\draw (5.6,-0.65) node{$\bar{t}$};

\draw[red]  (20,5.6) parabola (-6,0);

\draw[red]  (20,-5.6) parabola (-6,0);

\end{tikzpicture}
\caption{{\small Scheme of the proof of \eqref{bound}. The anisotropic paraboloid $\P^+_n$ (in red), centered in $P_n=(\,y,\, s-(\omega_n/C_1)^{2-\bar{p}}(C_2\rho_n)^{\bar{p}})$, evolves in a time $(\omega_n/C_1)^{2-\bar{p}}(C_2\rho_n)^{\bar{p}}$ to cover $\Q_{n+1}$.}}
 \label{FigA}
\end{figure}

\noindent To show that $\P_n^+\subset\Q_n\subset\Omega_T$, we control the space variables. From the upper bound in \eqref{doublebound} and of the paraboloid, we infer
\begin{equation*}
|x_i-y_i|^{p_i}<  \bigg(\frac{M_n-u(P_n)}{C_1}\bigg)^{p_i-2}\rho_n^{\bar{p}}\bigg(\frac{\omega_n}{C_1}\bigg)^{2-\bar{p}}\leq  \bigg( \frac{\omega_n}{C_1} \bigg)^{p_i-\bar{p}} \rho_n^{\bar{p}}
% \bigg(\frac{M_n-u(P_n)}{C_1}\bigg)^{p_i-2} \bigg(\frac{R}{A^n}\bigg)^{\bar{p}} \bigg(\frac{\omega_o}{C_1} \bigg)^{2-\bar{p}},
\end{equation*}
for all $x\in\pi_x(\P_n^+)$, being $\pi_x$ the projection on the space variables. This furnishes the desired inclusion. \\
% , still thanks to condition $\omega_n> M_n-u(P_n)$:
% \[|x_i-y_i|^{p_i}<  \bigg(\frac{M_n-u(P_n)}{C_1}\bigg)^{p_i-2}\rho_n^{\bar{p}}\bigg(\frac{\omega_n}{C_1}\bigg)^{2-\bar{p}} < \rho_n^{\bar{p}}\bigg(\frac{\omega_n}{C_1}\bigg)^{p_i-\bar{p}}.\]
\noindent
Now we show that, after a certain time $\bar{t}$, the whole cylinder $\Q_{n+1}$ is contained in the paraboloid $ \mathcal{P}^+_n$; see Figure \ref{FigA} for a representation. For times $t >s-(\omega_n/C_1)^{2-\bar{p}}(C_2\rho_n)^{\bar{p}}$, we denote by $\P^+_n(t)$ the time-section of $\P^+_n$ at time $t$:
\[ \P^+_n(t)= \bigg{\{} x \in \R^N: \, \,    |x_i-y_i|^{p_i}< C_2^{-\bar{p}} [(M_n-u(P_n))/C_1]^{p_i-2}(t-s+ (\omega_n/C_1)^{2-\bar{p}}(C_2\rho_n)^{\bar{p}})  \bigg{\}}.\]
\noindent Let us set \begin{equation*} \bar{t}=s-(\omega_{n+1}/C_1)^{2-\bar{p}}(C_2 \rho_{n+1})^{\bar{p}},\end{equation*}
% Indeed, for some $0<\rho<\rho_n$ this happens: we can determine $\rho$ by 
% \[(C_2\rho_{n+1})^{\bar{p}} (\omega_{n+1}/C_1)^{2-\bar{p}}= \rho^{\bar{p}} [(M_n-u(P_n))/C_1]^{2-\bar{p}}, \quad \Rightarrow \quad \rho=  \rho_{n+1} \bigg(\frac{M_n-u(P_n)}{\omega_n} \bigg)^{\frac{\bar{p}-2}{\bar{p}}}\leq \rho_n.\]
and let us prove that at time $\bar{t}$ we have the inclusion $\pi_{x}(\Q_{n+1})\subset \P^+_n(\bar{t})$. This reduces to show that
\[\rho_{n+1}^{\bar{p}} (\omega_{n+1}/C_1)^{p_i-\bar{p}} \leq (A^{\bar{p}}-1)  [(M_n-u(P_n)/C_1)]^{p_i-2} \rho_{n+1}^{\bar{p}} 
 (\omega_{n+1}/C_1)^{2-\bar{p}},\]
 that is,
\[\omega_{n+1}^{p_i-2} \leq (A^{\bar{p}}-1) (M_n-u(P_n))^{p_i-2}.\]
According to \eqref{doublebound}, this inequality is verified when $4^{p_N-2}<A^{\bar{p}}-1$, as for instance setting $A= 4^{p_N}$. \vskip0.1cm \noindent Hence, by the Harnack inequality \eqref{estimate-paraboloid} and \eqref{doublebound}, we can estimate the infimum of $M_{n}-u$ in $\Q_{n+1}$ as
\begin{equation} \label{bound}
 \inf_{\Q_{n+1}} (M_{n}-u) \ge \inf_{\P^+_n(\bar{t})} (M_n-u) \ge  \frac{M_{n}-u(P_{n})}{C_3} \ge\, \,  \omega_{n+1}/(4C_3),
\end{equation}\noindent again referring to Figure \ref{FigA}. Thus
\[
M_{n} \ge \sup_{\Q_{n+1}} u+ \omega_{n+1}/(4C_3).  
\]
Adding $-\inf_{\Q_n} u \ge - \inf_{\Q_{n+1}}u$ to both sides, besides using $\osc_{\Q_{n+1}}u >\omega_{n+1}$, we get 
\[\omega_n \ge M_n - \inf_{\Q_n} u \ge  \sup_{\Q_{n+1}} u+ \omega_{n+1}/(4C_3)  - \inf_{\Q_{n+1}} u= \osc_{\Q_{n+1}}u+\omega_{n+1}/(4C_3)> \bigg(1+\frac{1}{4C_3}  \bigg) \omega_{n+1}\, .
\]This leads to a contradiction by definition of $\delta$, since
\[
\omega_n > \bigg(1+\frac{1}{4C_3}  \bigg) \delta \omega_{n}=  \bigg(\frac{4C_3}{1+4C_3}\bigg) \bigg(1+\frac{1}{4C_3} \bigg) \omega_n = \omega_n.
\] 

\end{proof}
\noindent 
{\small STEP 4-{\it Conclusion of the proof of Theorem \ref{HC}.}}
\vskip0.2cm 

\noindent If we consider a point $(x,t) \in (y,s) +\Q_R^-(\omega_o/C_1,C_2)$, let $n \in \mathbb{N}$ be the last number such that we have $(x,t)\in \Q_n$, so that $(x,t) \not\in \Q_{n+1}$. From the first condition and \eqref{control} we have 

\[|u(x,t)-u(y,s)| \leq \osc_{\Q_n} u\leq \delta^n \omega_o.\]
The rest of the job is standard and consists in determining from condition $(x,t) \not\in \Q_{n+1}$ an upper bound for $\delta^n$. For the sake of simplicity, we just show the case $x \not\in y+\K_{\rho_{n+1}}$.

% To this aim, we distinguish two cases.
% \vskip0.1cm
% \noindent \textit{Case 1:} $x \not\in y+\K_{\rho_{n+1}}$.\newline
\noindent 
Let $\beta>0$ be such that $\delta^{{(\bar{p}-2)}/{\bar{p}}}/A=\delta^{\beta}$. By assumption, there is an index $i \in \{1,\dots, N\}$ such that
\begin{equation*}
%= \bigg(\frac{\delta^{\frac{\bar{p}-2}{\bar{p}}}}{A}\bigg)^{n{\frac{\bar{p}}{p_i}}}&  \bigg(\frac{\delta^{\frac{\bar{p}-2}{\bar{p}}}}{A}\bigg)^{\frac{\bar{p}}{p_i}} R^{\frac{\bar{p}}{p_i}} (\delta^n \omega_o/C_1)^{\frac{p_i-\bar{p}}{p_i}}\\
% &
%&|x_i-y_i|> \rho_{n+1}^{{\bar{p}}/{p_i}} (\omega_{n+1}/C_1)^{{(p_i-\bar{p})}/{p_i}}
%\ge \gamma(\delta,A) (\delta^{n})^{[{\bar{p}(\beta-1)+p_i}]/{p_i}}R^{{\bar{p}}/{p_i}} (\delta^n \omega_o/C_1)^{{(p_i-\bar{p})}/{p_i}}, \\
|x_i-y_i|^{p_i}> \rho_{n+1}^{{\bar{p}}} (\omega_{n+1}/C_1)^{{p_i-\bar{p}}} \geq \gamma(A)(\delta^n)^{[{\bar{p}(\beta-1)+p_i}]}R^{\bar{p}}(\omega_o/C_1)^{p_i-\bar{p}},
\end{equation*}
\noindent that gives us, for $\chi_i=\bar{p}/(\bar{p}(\beta-1)+p_i)$, the following estimate of $\delta^n$:
 \begin{equation*} \begin{aligned}
 \delta^n \leq&
 \gamma  \bigg( \frac{ |x_i-y_i|^{{p_i}/{\bar{p}}}  (\omega_o/C_1)^{{(\bar{p}-p_i)}/{\bar{p}}}}{R} \bigg)^{{\bar{p}}/[{{\bar{p}(\beta-1)+p_i}}]} \\
 &\leq \gamma \bigg(\frac{\sum_{i=1}^N |x_i-y_i|^{{p_i}/{\bar{p}}}\omega_o^{{(\bar{p}-p_i)}/{\bar{p}}}+ |t-s|^{{1}/{\bar{p}}}\omega_o^{{(\bar{p}-2)}/{\bar{p}}}}{{\bf{p}}\text{-dist}(K,\partial \Lambda) } \bigg)^{\chi_i}.\end{aligned} \end{equation*}
 \noindent
 From $A>4>\delta^{-1-2/\bar{p}}$ we infer $\beta>2$, whence $\chi_i\in(0,1)$. A similar estimate follows from the case where times are not contained, with $\chi_t= \bar{p}/(\bar{p}(\beta-1)+2)$. Therefore, recalling that $p_N>2$, we choose the H\"older exponent
\begin{equation} \label{alfa}
    \chi = \min \{\chi_i, \chi_t, \quad i=1,\dots,N \}=\frac{\bar{p}}{\bar{p}(\beta-1)+p_N}.
\end{equation}
 \end{proof}

\noindent 

\section{Liouville-type results}\label{LiouvilleSection}
\noindent In their origins, Liouville properties were discovered for harmonic functions. Indeed, for solutions to $\Delta u =0$ in $\R^N$, a one-sided bound on $u$ or the sublinear growth at infinity are suitable rigidity conditions.  
% The easiest known ones are, for harmonic functions 
% \begin{itemize}
%     \item[(i)] Any solution $u$ of $\Delta u=0$ in $\R^N$ which is bounded from below is constant;
%     \item[(ii)] Any solution $u$ of $\Delta u=0$ in $\R^N$ which grows sublinearly at infinity is constant.
% \end{itemize}
These two classical examples  follow respectively from an application of Harnack's inequality and from gradient estimates. Here we observe that gradient bounds of logarithmic type are unknown for solutions to the stationary counterpart of \eqref{EQ} and seem hard to obtain, chiefly because of the lack of homogeneity of the operator. On the other hand, for parabolic equations a one-side bound is not sufficient to imply that solutions are constant, as we remarked. This is still the case also for non-negative solutions to degenerate $p$-Laplacian equations (i.e., for $p>2$). Indeed, the one-parameter family of non-negative functions 
\[
\R \times \R  \ni (x,t)\rightarrow u(x,t;c)= c^{{1}/({p-2})} \bigg(\frac{p-2}{p-1}  \bigg)^{{(p-1)}/{(p-2)}} (1-x+ct)_+^{{(p-1)}/{(p-2)}}
\] is a family of non-negative, non-constant weak solutions to $u_t=\Delta_{p}u$ in $\R^2$. This naturally provides a counterexample also in case of equation \eqref{EQ} in one spatial dimension. Similarly, the anisotropic driving example we have in mind is
\[\R^N \times \R \ni (x,t) \rightarrow u(x,t;c)= \bigg(1-ct + \sum_{i=1}^N {(\alpha_i/p_i') |x_i|^{p_i'}}\bigg)_+,\]
for $\alpha_i >0$ such that $\sum_{i=1}^N |\alpha_i|^{p_i-1}\alpha_i=c$ and being $p_i'$ the H\"older conjugate of $p_i$ for each $i=1,\dots, N$. On the other hand, a full lower bound coupled with a specific upper bound at some time level ensures a Liouville property, as the following result uncovers.

\begin{theorem}\label{Liouville1}
Let $T\in \R$, $S_T=\R^N \times (-\infty, T)$, and $u$ be a solution to \eqref{EQ}-\eqref{pi} which is bounded below in $S_T$. Assume moreover that, for some $s<T$, one has 
\begin{equation}\label{rs-bound}
    \sup_{\R^N} u(\cdot,s)=M_s <\infty.
\end{equation} Then $u$ is constant in $S_s= \R^N \times (-\infty, s)$.
\end{theorem} 
\begin{corollary}
\label{LiouvilleCor}
Let $T\in \R$, $S_T=\R^N \times (-\infty, T)$, and $u$ be a solution to \eqref{EQ}-\eqref{pi}. If $u$ is bounded from above and below in $S_T$, then it is constant.
\end{corollary}

\begin{proof}[Proof of Theorem \ref{Liouville1}]
Let $u$ be a solution to \eqref{EQ} bounded from below in $S_T$. We define
\[
m:=\inf_{S_T} u.
\]
%For any point $(y,s)\in S_T$ such that $u(y,s)>m$, we set $\theta:=\frac{u(y,s)-m}{C_1}$ and consider the intrinsic backward $\bf{p}$-paraboloid $\P_\theta^-(y,s)$ as in Section 2.
We prove the following fact, which is interesting in its own:
\begin{equation}\label{perse}
    \lim_{t\rightarrow -\infty} u(x,t)= \inf_{S_T} u \quad \mbox{for any} \quad x \in \R^N.
\end{equation}\noindent To this aim, fix any $x\in\R^N$ and $\varepsilon>0$. Notice that there exists a point $(y_\varepsilon, s_\varepsilon) \in S_T$ such that $u(y_{\varepsilon},s_{\varepsilon})-m\leq \varepsilon/C_3$. Set $\theta_{\varepsilon}= (u(y_{\varepsilon},s_{\varepsilon})-m)/C_1$. Exploiting \eqref{estimate-paraboloid} for the solution $u-m$, we have
\begin{equation}\label{infimum}
m \leq u(y,s) \leq m+\varepsilon \quad \text{for all} \quad (y,s) \in \P_{\theta_\varepsilon}^-(y_{\varepsilon},s_{\varepsilon}).
\end{equation}
Consider the half line $R:=\{x\}\times(-\infty,T)$. Observe that
\[
R \cap \P_{\theta_\varepsilon}^-(y_{\varepsilon},s_{\varepsilon}) = \{x\} \times (-\infty,t_{\varepsilon,x}), \quad \mbox{being} \quad t_{\varepsilon,x} := s_\varepsilon-C_2^{\bar{p}}(2-|x_i-y_{\varepsilon,i}|)^{p_i}\theta^{2-p_i}.
\]
According to \eqref{infimum}, this shows that
\begin{equation*}
m \leq u(x,s) \leq m+\varepsilon \quad \text{for all} \quad s < t_{\varepsilon,x}.
\end{equation*}
Accordingly, \eqref{perse} is proved, by arbitrariness of $x$ and $\varepsilon$. 
% Recalling the definition of $t_{\varepsilon,x}$, it is readily seen that \eqref{perse} is actually uniform in $x$ on compact subsets of $\R^N$. \newline
A similar argument shows that
\begin{equation}\label{perse2}
    \sup_{S_T} u <\infty \quad \Rightarrow \quad \lim_{t\rightarrow -\infty} u(x,t) = \sup_{S_T} u \quad \forall x \in \R^N.
\end{equation}
Eventually this implies that any $u$ solution to \eqref{EQ} which is bounded from both above and below in the whole $S_T$ is necessarily constant. Indeed, by \eqref{perse} and \eqref{perse2} we have $\sup_{S_T} u= \inf_{S_T} u$. This argument proves Corollary \ref{LiouvilleCor}.
\vskip0.1cm \noindent 
In order to conclude the proof of Theorem \ref{Liouville1}, we use the assumption that there exists $\bar{s} \in (-\infty,\, T)$ such that $u(\cdot, \bar{s})$ is bounded from above in the whole $\R^N$ by a suitable $M_s\in \R$. Indeed, letting $\theta_x=(u(x,\bar{s})-m)/C_1$ for any $x\in\R^N$ and using the intrinsic backward Harnack inequality for $u-m$ again, we get the uniform bound 
\[
u(y,s) \leq C_3 u(x,\bar{s}) \leq C_3 M_{\bar{s}}, \quad \text{for all} \quad x \in \R^N \quad \mbox{and} \quad (y,s) \in \P_{\theta_{x}}^-(x,\bar{s}).
\] Reasoning as above, with $\P_{\theta_x}^-(x,\bar{s})$ instead of $\P_{\theta_{\varepsilon}}^- (y_{\varepsilon},\, s_{\varepsilon})$, besides recalling that $u$ bounded from both above and below in $\P_{\theta_x}^-(x,\bar{s})$ uniformly in $x \in \R^N$, we conclude that $u$ is constant in $S_{\bar{s}}$.
\end{proof} 

\noindent As a general principle, the bigger the set where the equation is solved the stronger the rigidity: for solutions of \eqref{EQ} in $\R^N \times \R$, it suffices to check their asymptotic (in time) two-side boundedness at a single point $y \in \R^N$ to infer that they are constant, as shown by the next theorem. 
\begin{theorem}\label{Liouville2}
Let $u$ be a local weak solution to \eqref{EQ}-\eqref{pi} in $\R^N \times \R$ which is bounded from below. If, in addition, there exists $y \in \R^N$ and a sequence $\{ s_n\} \subset \R$, $s_n\to+\infty$, such that $\{u(y, s_n)\}$ is bounded, then $u$ is constant.
\end{theorem}\noindent 
\begin{remark}
We explicitly point out the following straightforward consequence of Theorem \ref{Liouville2}. Let $u$ be a local weak solution to \eqref{EQ} in $\R^N \times \R$ which is bounded from below. Suppose that, for some $y \in \R^N$, one has
\begin{equation}\label{Liouville2HP}
    \liminf_{t\rightarrow +\infty} u(y,t)=\alpha \in \R.
\end{equation}\noindent Then $u$ is constant.
\end{remark}

\begin{proof}[Proof of Theorem \ref{Liouville2}]
Let $m:=\inf u$ and consider $\tilde{u}:=u+m+C_1$, which is a solution to \eqref{EQ}. By assumption, there exist $M \in \R$ and $\{s_n\}\subset\R$ such that $s_n \to +\infty$ and
\[
\tilde{u}(y,s_n)<M \quad  \quad \forall n \in \N.
\]
Let us fix arbitrarily $\bar{s}\in \R$ and let $\bar{n} \in \N$ be big enough such that $s_n >\bar{s}$ for all $n \ge \bar{n}$. Then, for all $n\ge\bar{n}$, we set $\theta_n:=\tilde{u}(y,s_n)/C_1$ and define a sequence of radii $\{\rho_n\}$ through
\[
s_n-\theta_n^{2-\bar{p}} (C_2 \rho_n)^{\bar{p}} =\bar{s}, \quad \quad \mbox{that is,} \quad \quad \rho_n= [\theta_n^{\bar{p}-2} (s_n-\bar{s})]^{1/\bar{p}}/C_2.
\] We want to apply the Harnack inequality to deduce an upper bound for $\tilde{u}(\cdot,\bar{s})$ in the whole $\R^N$; so we need to check that the intrinsic anisotropic cubes $\K_{\rho_n}(\theta_n)$ expand as $s_n\to+\infty$. An explicit computation yields
\[
\K_{\rho_n}(\theta_n)= \prod_{i=1}^N \bigg{\{}|x_i|<\theta_n^{{(p_i-2)}/{p_i}}  \left(\frac{s_n-\bar{s}}{C_2^{\bar{p}}}\right)^{{1}/{p_i}} \bigg{\}}\quad \xrightarrow[n \to \infty]{} \quad \R^N,
\] since $1\le\theta_n\le M/C_1$ and $\{s_n\}$ diverges. By the intrinsic Harnack inequality \eqref{Harnack} we have
\[
\sup_{y+\K_{\rho_n}(\theta_n)} \tilde{u} \bigg(\, \,  \cdot\, \, ,\,  s_n-\theta_n^{2-\bar{p}} (C_2\rho_n)^{\bar{p}} \bigg)\leq C_3 \,\tilde{u}(y,s_n) \leq C_3 M \quad \forall n\ge \bar{n}.
\] Thus, recalling the definition of $\{\rho_n\}$, we get the uniform estimate 
\[
\sup_{y+\K_{\rho_n}(\theta_n)} \tilde{u}(\cdot, \bar{s} ) \leq C_3 M, \qquad \forall n\ge \bar{n},
\]
whence, letting $n\to\infty$,
\[
\sup_{\R^N} \tilde{u}(\cdot, \bar{s}) \leq C_3 M.
\]
Now we can apply Theorem \ref{Liouville1} in $(-\infty, \bar{s})$ and conclude by the arbitrariness of $\bar{s}\in\R$.
\end{proof}

\noindent Finally, we show that the oscillation estimates \eqref{control} constitute a Liouville property for ancient solutions. This allows us to get rid of the range of $p_i$s of finite speed of propagation \eqref{pi}, at the price of assuming a suitable decay of the local oscillation.

\begin{theorem} \label{Cutilisci}
Let $u$ be a bounded function in $S_T$. Let $\omega_o>0$, $\delta \in (0,1)$, $2<p_1 \leq \ldots \leq p_N< \infty$ and $c_1,c_2,c_4>1$ be fixed parameters. For any $(\bar{x}, \bar{t}) \in S_T$ and $R_o>0$, define a sequence of backward shrinking cylinders $\Q_{n+1} \subset \Q_{n}$ as  
\begin{equation}\label{cilindrotti}
   \Q_n = \Q_{n}(R_o) = (\bar{x}, \bar{t})+ \Q_{\rho_n}^-(\theta_n, c_2) = \prod_{i=1}^N \bigg{\{}|x_i-\bar{x}_i|<\theta_n^{{(p_i-\bar{p})}/{p_i}}\rho_n^{{\bar{p}}/{p_i}}  \bigg{\}} \times \bigg(\bar{t}-\theta_n^{2-\bar{p}} (c_2\rho_n)^{\bar{p}} ,\, \bar{t}\, \bigg],
\end{equation} being
\[
\theta_n= \delta^n \omega_o/c_1, \qquad \rho_n=\varepsilon^n R_o, \qquad \varepsilon=\delta^{\frac{\bar{p}-2}{\bar{p}}}/c_4.
\]
\noindent If $u$ satisfies, for all $(\bar{x},\bar{t}) \in S_T$ and $R_o>0$, the decay
\begin{equation}\label{oscilla}
    \osc_{\Q_{n+1}} u \leq \delta \osc_{\Q_{n}} u, \qquad n\in \N \cup\{0\}, 
\end{equation} 
\noindent then $u$ is constant in $S_T$.
\end{theorem}

\begin{proof} The proof is an adaptation of an early idea already present in \cite{Glagoleva1} (see also \cite{Landis}). Arguing by contradiction, assume that $A,B \in S_T$ are two points such that $u(A) \ne u(B)$. %and call $Z_T= \{T\} \times \R^N$.
Suppose, without loss of generality, $T=0$ and define
\[d= \max \{\dist(A,0), \, \dist(B,0) \}.\]
%\[d= \max \{\dist(A, Z_T),\,  \dist(B, Z_T),\, \dist(A,B), \, \dist(A,0), \, \dist(B,0) \}.\]
Choose a radius $\tilde{R}_o>0$ big enough to enclose $A$ and $B$ inside an intrinsic backward cylinder $\tilde{\Q}_0:=\Q_0(\tilde{R}_o)$, %$\tilde{Q}_{R_0}^-$
so that $\tilde{R}_o$ satisfies 
\begin{equation*}
    \begin{cases}
  \theta_0^{p_i-\bar{p}} \tilde{R}_o^{{\bar{p}}} >d^{p_i},\quad i =1,..,N,\\
\theta_0^{2-\bar{p}} (c_2 \tilde{R}_o)^{\bar{p}} >d.
    \end{cases}% \quad  \text{with} \quad \delta = 4C_3/(1+4C_3),
\end{equation*} %where $C_1,C_2,C_3$ are the constants of Theorem \ref{Harnack-Inequality}.
\noindent Now set $R_o:=c_4 \tilde{R}_o\delta^{\frac{2-p_N}{\bar{p}}}$, observe that $\tilde{\Q}_0\subset\Q_1(R_o)$, and fix $\tilde{\Q}_1:=\Q_0(R_o)$.
%be defined as
%\[
%\tilde{{\Q}}_1= \prod_{i=1}^N \bigg{\{}|x_i|<[\delta^{-1}(\omega_o/c_1)]^{{(p_i-\bar{p})}/{p_i}} (\tilde{R_o}/\varepsilon)^{{\bar{p}}/{p_i}}  %\bigg{\}} \times \bigg(-[\delta^{-1}(\omega_o/c_1)  ]^{2-\bar{p}} [(c_2 \tilde{R_o})/\varepsilon ]^{\bar{p}}\,  ,\, 0\, \, \bigg],
%\] 
%whence $\tilde{\Q}_0 \subset \tilde{\Q}_1$.
Then the decay \eqref{oscilla} implies
\[\osc_{\tilde{\Q}_0} u \le \osc_{\Q_1(R_o)} u \le \delta \osc_{\Q_0(R_o)} u = \delta \osc_{\tilde{\Q}_1} u.\]
Proceeding inductively, we construct $\tilde{\Q}_{n+1}$ by choosing a new $R_o$ such that $\tilde{\Q}_n \subset \Q_1(R_o) \subset \Q_0(R_o) =: \tilde{\Q}_{n+1}$. By construction,
\[ |u(A)-u(B)| \le \osc_{\tilde{\Q}_0} u \le \delta^n \osc_{\tilde{\Q}_{n}} u  \quad  \quad \forall n \in \N.\]
Finally, the boundedness of $u$ leads to a contradiction: indeed, for all $n \in \N \cup \{0\}$ we have
\[|u(A)-u(B)| \le \delta^n \osc_{\tilde{\Q}_{n}} u \le 2 \delta^n \|u\|_{\infty, S_T},\]
forcing $u(A)=u(B)$.
\end{proof}

%\begin{proof}[Five-lines-proof]
%By iterating \eqref{oscilla} and the boundedness of $u$, we get 
%\[\osc_{\Q_n} \leq \delta^n \osc_{\Q_0} u \leq 2 \delta^n \|u\|_{\infty,S_T}.\]
%Passing to the limit $R\rightarrow \infty$ first, we obtain
%\[ \osc_{S_T} u \leq 2 \delta^n \|u\|_{\infty,S_T},\]
%and then passing to the limit for big $n$ finishes the proof. 
%\end{proof}
\noindent As a consequence of Proposition \ref{birra}, we obtain again Corollary \ref{LiouvilleCor}. Indeed, when equation \eqref{EQ} is solved in $S_T$ the length $R$ in Proposition \ref{birra} can be taken arbitrarily large. Nevertheless, we decided to formulate  Theorem \ref{Cutilisci} without the assumption that $u$ is a solution of any equation. Indeed, Theorem \ref{Cutilisci} is finer: its general principle goes far beyond equation \eqref{EQ} and is a key argument to prove rigidity results for a very general class of equations (see, e.g., \cite[Prop. 18.4]{DBGV-mono} or, for instance, \cite{Liao2} for an application to systems). Its importance shows up when a Harnack inequality ceases to hold true.

\section{Time-extrinsic Harnack inequality} \label{WHSection}
\noindent In this section we show how it is possible to free the Harnack inequality from its intrinsic geometry in time. More specifically, we give a formulation of the Harnack inequality allowing the solution to be evaluated at any time level, independently of the anisotropic geometry, provided there is enough room for the anisotropic evolution inside $\Omega_T$. Unlike the isotropic case, here it looks harder to get rid of the intrinsic geometry along the space variables. The proof of the next theorem exploits a comparison with the abstract Barenblatt solution $\B$ of Theorem \ref{Barenblatt} to control the positivity.

\begin{theorem}
\label{WeakHarnack} Let $u \ge 0$ be a local weak solution to \eqref{EQ} in $\Omega_T$, and assume \eqref{pi}. Then there exist $ \tilde{\eta} >0$ and $\gamma>1$, depending only on $N$ and $p_i$s, such that for all $(x_o,t_o)\in \Omega_T$ and $\rho, \tilde{\theta}>0$ fulfilling the condition
\begin{equation}\label{domain}
(x_o, t_o+ \tilde{\theta})+\Q_{C_3 \rho}(u(x_o,t_o)/C_1,C_2) \subset \Omega_T
\end{equation} we have
\begin{equation}\label{WH}
    u(x_o,t_o)\leq \gamma \bigg{\{} \bigg(\frac{\rho^{\bar{p}}}{\tilde{\theta}}  \bigg)^{{1}/{(\bar{p}-2)}}+ \bigg( \frac{\tilde{\theta}}{\rho^{\bar{p}}} \bigg)^{N/\bar{p}} \bigg[\inf_{x_o+K_{\tilde{\eta}\rho}(\tilde{\eta} u(x_o,t_o)/C_1)} u( \cdot,\, t_o+\tilde{\theta})  \bigg]^{\lambda/\bar{p}}\bigg{\}},
\end{equation} \noindent where $C_1,C_3>1$ come from Theorem \ref{Harnack-Inequality} while $\lambda,\tilde{\eta}>0$ stem from Theorem \ref{Barenblatt}. \end{theorem}
\begin{proof}
Let $\rho, \tilde{\theta} >0$ be such that \eqref{domain} holds true. Set
\begin{equation} \label{t*}
t^*:= \bigg( \frac{C_1}{u(x_o,t_o)} \bigg)^{\bar{p}-2} (C_2 {\rho})^{\bar{p}}.
\end{equation} We can suppose $t^*<\tilde{\theta}/2$; otherwise we get $u(x_o,t_o) \leq \gamma ({\rho}^{\bar{p}}/\tilde{\theta})^{1/(\bar{p}-2)}$ for a suitable $\gamma= \gamma (C_1, C_2,\bar{p})$, and \eqref{WH} is valid. Observe that $t^*<\tilde{\theta}/2$ and \eqref{domain} imply
\[
t_0+\bigg( \frac{C_1}{u(x_o,t_o)}\bigg)^{\bar{p}-2} (C_2 {\rho})^{\bar{p}} <t_0+ \tilde{\theta}/2 < T \quad \mbox{and} \quad x_o + \K_{C_3 {\rho}}(u(x_o,t_o)/C_1) \subset \Omega.
\] Hence the forward Harnack inequality \eqref{Harnack} furnishes
\[
u(x_o,t_o) \leq C_3 u(x,\, t_o+t^*) \quad \quad \forall \, x \in x_o + \K_{{\rho}}(u(x_o,t_o)/C_1).
\] This initial value can be considered for a comparison with the Barenblatt solution  $\B_{\sigma}(x-x_o,t-s)$ centered at $(x_o,s)$, being $s,\sigma>0$ to be chosen such that $\B_{\sigma}(x-x_o,t_o+t^*-s)$ lies below $u$ in $x_0+\K_{\rho}(u(x_o,t_o)/C_1)$. These requirements can be written as
\begin{equation}\label{supportami}
\begin{cases}
\supp{\B_{\sigma}}(\cdot-x_o, t_o+t^*-s) \subseteq x_o+ \K_{\rho}(u(x_o,t_o)/C_1),\\
\|\B_{\sigma}(\cdot-x_o, t_o+t^*-s)\|_\infty \leq u(x_o,t_o)/C_3.
\end{cases}
\end{equation}
According to Theorem \ref{Barenblatt}, conditions in \eqref{supportami} are fulfilled as long as
\begin{equation}\label{supportami2}
\begin{cases}
 \sigma^{(p_i-2)/p_i} (t_o+t^*-s)^{\alpha_i} \leq \rho^{\bar{p}/p_i} (u(x_o,t_o)/C_1)^{(p_i-\bar{p})/p_i},\\
 \sigma (t_o+t^*-s)^{-\alpha} \leq u(x_o,t_o)/C_3.
\end{cases}
\end{equation}
\noindent Inequalities in \eqref{supportami2} are in turn ensured by choosing
\[
\sigma= (t_o+t^*-s)^{N/\lambda} u(x_o,t_o)/C_3 \quad \mbox{and} \quad s= t_o+t^*- \bigg( \frac{\rho^{\bar{p}}}{u(x_o,t_o)^{\bar{p}-2}} \bigg) \gamma_1,\] 
where $\gamma_1=\min \{(C_3^{p_i-2})/(C_1^{p_i-\bar{p}})\, |\, i=1,\dots,N\}$. Therefore the comparison principle, applied at the time $t_o+\tilde{\theta}>t_o+t^*$, gives
\begin{equation}\label{comparison} \begin{aligned}
u(x, t_o+\tilde{\theta}) &\ge \tilde{\eta} \sigma |t_o+t^*-(t_o+\tilde{\theta})|^{-\alpha} = \tilde{\eta} \bigg( \frac{u(x_o,t_o)}{C_3} \bigg) (t_o+t^*-s)^{N/\lambda} (\tilde{\theta}-t^*)^{-N/\lambda}\\
&\ge \tilde{\eta} \bigg( \frac{u(x_o,t_o)}{C_3} \bigg) \bigg(  \frac{\gamma_1  \rho^{\bar{p}}}{u(x_o,t_o)^{\bar{p}-2}} \bigg)^{N/\lambda} \tilde{\theta}^{-N/\lambda} \ge \gamma u(x_o,t_o)^{\bar{p}/\lambda} \bigg( \frac{\rho^{\bar{p}}}{\tilde{\theta}} \bigg)^{N/\lambda},
\end{aligned} \end{equation} with $\gamma= \gamma(\gamma_1,\tilde{\eta})$, for every $x$ in the set of positivity
\begin{equation*}
\begin{aligned}
\mathcal{P}_{t_o+\tilde{\theta}-s}(x_o)\supseteq \mathcal{P}_{t_o+t^*-s}(x_o) &= \prod_{i=1}^N\{|x_i-x_{o,i}|\leq\tilde{\eta} \rho^{\bar{p}/p_i} (u(x_o,t_o)/C_1)^{(p_i-\bar{p})/p_i} \}\\
&=x_o+ \K_{\tilde{\eta}\rho}(\tilde{\eta}u(x_o,t_o)/C_1),
\end{aligned}
\end{equation*}
with a constant $\tilde{\eta}$ depending only on the data $N, p_i$. Taking the infimum in the estimate \eqref{comparison} on the set $x_o+ \K_{\tilde{\eta}\rho}(\tilde{\eta}u(x_o,t_o)/C_1)$ concludes the proof.
\end{proof}

\begin{remark}
In Theorem \ref{WeakHarnack} the lower bound $u(x_o,t_o)>0$ is not required; moreover, $\tilde{\theta}>0$ is arbitrarily chosen between those numbers that preserve the inclusion \eqref{domain}. When the equation is solved in $\R^{N+1}$, the proof furnishes inequality \eqref{WH} without the first term on the right. 
% Henceforth, when the equation is solved in $\R^{N+1}$, the proof furnishes inequality \eqref{WH} without the first term on the right, and one can infer similar Liouville properties as Theorem \ref{Liouville2}.
\vskip0.1cm \noindent Actually, Theorems \ref{Harnack-Inequality} and \ref{WeakHarnack} are equivalent for small radii. We can easily show that Theorem \ref{WeakHarnack} implies Theorem \ref{Harnack-Inequality} by a simple choice of $\tilde{\theta}$. For instance, let us pick
\[ \tilde{\theta}= (2\gamma)^{\bar{p}-2}\rho^{\bar{p}}u(x_o,t_o)^{2-\bar{p}},
\] and suppose that $(x_o,t_o+\tilde{\theta}) + \Q_{C_3\rho} (u(x_o,t_o)/C_1)\subset \Omega_T$. Then inequality \eqref{WH} leads to
\[
u(x_o,t_o) \leq \gamma \bigg{\{} \frac{u(x_o,t_o)}{2\gamma}+ \bigg( \frac{2\gamma}{u(x_o,t_o)} \bigg)^{{N(\bar{p}-2)}/{\bar{p}}}\bigg[ \inf_{x_o+ \K_{\tilde{\eta} \rho}(\tilde{\eta} u(x_o,t_o)/C_1)} u(\, \cdot\, , \, t_o+ \bigg( \frac{u(x_o,t_o)}{2\gamma} \bigg)^{2-\bar{p}} \rho^{\bar{p}}) \bigg]^{{\lambda}/{\bar{p}}}  \bigg{\}},
\]
whence
\[
u(x_o,t_o) \leq \tilde{C_3} \inf_{x_o+ \K_{\tilde{\rho}}(M)} u(\cdot, \, t_o+ \tilde{C_2} M^{2-\bar{p}} \tilde{\rho}^{\bar{p}}), \quad M= u(x_o,t_o)/\tilde{C}_1, \] for all $\tilde{\rho} \leq \tilde{\eta} \rho$ and with positive constants 
\[\tilde{C_1}= C_1/\tilde{\eta}, \quad  \tilde{C_2}=\tilde{\eta}^{-2}(2\gamma/\tilde{C}_1)^{\bar{p}-2}, \quad \tilde{C_3}= 2 \gamma.\] \end{remark}

\section*{Acknowledgements}
\noindent 
We are grateful to S.A. Marano and V. Vespri for encouraging us toward this project. We wish to thank professor S. Mosconi for his precious suggestions and E. Macca for a numerical insight about Barenblatt-type solutions. Moreover, we are indebted with E. Henriques for pointing out an early mistake about H\"older continuity of solutions. Finally, S. Ciani is supported by the department of Mathematics of Technical University of Darmstadt, and U. Guarnotta is supported by: (i) PRIN 2017 `Nonlinear Differential 
Problems via Variational, Topological and Set-valued Methods' (Grant 
No. 2017AYM8XW) of MIUR; (ii) GNAMPA-INdAM Project 
CUP$\underline{\phantom{x}}$E55F22000270001; (iii) grant `PIACERI 
20-22 Linea 3' of the University of Catania.

\end{document}